\documentclass{amsproc}
\usepackage{amsmath,amsthm, amscd, amssymb, amsfonts}

\usepackage[all]{xy}

\usepackage{hhline}

\usepackage[active]{srcltx}
\usepackage{color}

\renewcommand{\_}[1]{_{\left( #1 \right)}}
\renewcommand{\^}[1]{^{\left( #1 \right)}}

\newcommand\sk{\mathbb S_4}
\def\unon{\left\{ 1,\dots ,\theta \right\}}

\newcommand{\fd}{finite-dimensional}

\newcommand{\sv}{\mathcal{SV}}

\newcommand{\ms}{\mathcal{SM}}
\newcommand{\supre}{{}_H\hspace{-1pt}\ms }
\newcommand{\supa}{{}_A\hspace{-1pt}\ms }
\newcommand{\supc}{{}^C\hspace{-1pt}\ms }
\newcommand{\supcr}{\ms^C }

\newcommand{\supaa}{\supa\hspace{-1pt}_A}

\newcommand{\mas}{{}_{A^\sigma}\hspace{-1pt}\mathcal M }

\newcommand{\cocs}{{}^{C^\sigma}\hspace{-3pt}\mathcal M }

\newcommand{\supcomh}{{}^H\hspace{-2pt}\ms }

\newcommand{\suphbim}{{}^H_H\hspace{-2pt}\ms^H_H }

\newcommand{\supertau}{\tau_{\text{super}}}

\newcommand{\co}{\operatorname{co} }

\def\de{\Delta}
\newcommand\ord{\operatorname{ord}}

\newcommand{\otv}{{1\le i \le \theta}}
\newcommand{\otvz}{{1\le i, j \le \theta}}
\newcommand{\otrvz}{{1\le i \neq  j \le \theta}}

\newcommand{\e}{\mathbf e}
\newcommand{\f}{\mathbf f}

\newcommand{\ub}{\mathbf u}
\newcommand{\vb}{\mathbf v}
\newcommand{\wb}{\mathbf w}

\newcommand{\zb}{\mathbf z}

\newcommand{\si}{s_{F,i}}
\newcommand{\sE}{s_{E,i}}

\newcommand{\Ss}{{\mathcal S}}

\newcommand{\hs}{H^{\sigma}}

\def\cX{\mathcal{X}}

\newcommand{\J}{{\mathcal J}}
\newcommand{\R}{{\mathcal R}}

\newcommand{\incl}{{\mathfrak i}}

\def\zt{\Z^{\theta}}

\newcommand\ot{\otimes}

\newcommand{\ku}{\Bbbk}
\newcommand{\kk}{\mathbb C}

\newcommand{\Z}{{\mathbb Z}}
\newcommand{\N}{{\mathbb N}}

\newcommand{\G}{\mathbb{G}}

\def\cG{\mathcal{G}}

\newcommand{\toba}{{\mathcal B}}
\newcommand{\T}{{\mathcal T}}
\newcommand{\Hc}{{\mathcal H}}

\newcommand{\ydk}{{}^K_K\mathcal{YD}}

\newcommand{\ydh}{{}^H_H\mathcal{YD}}
\newcommand{\ydsh}{{}^H_H\mathcal{YDS}}
\newcommand{\ydsho}{{}^{H_0}_{H_0}\mathcal{YDS}}

\newcommand{\ydg}{{}^{\ku\Gamma}_{\ku\Gamma}\mathcal{YD}}

\newcommand{\End}{\operatorname{End}}
\newcommand{\Aut}{\operatorname{Aut}}

\newcommand\gr{\operatorname{gr}}

\def\cR{\mathcal{R}}
\def\cC{\mathcal{C}}

\newcommand\sgn{\operatorname{sgn}}
\newcommand\ad{\operatorname{ad}}
\newcommand\coin{\operatorname{co}}
\newcommand\Hom{\operatorname{Hom}}

\newcommand{\ydsg}{{}^{\ku\Gamma}_{\ku\Gamma}\mathcal{YDS}}

\def\cW{\mathcal{W}}
\numberwithin{equation}{section}\theoremstyle{plain}

\newtheorem{theorem}{Theorem}[section]
\newtheorem{lema}[theorem]{Lemma}

\newtheorem{prop}[theorem]{Proposition}

\theoremstyle{definition}
\newtheorem{definition}[theorem]{Definition}
\newtheorem{exa}[theorem]{Example}

\theoremstyle{remark}
\newtheorem{obs}[theorem]{Remark}

\newcommand\id{\operatorname{id}}

\def\pf{\begin{proof}}
\def\epf{\end{proof}}

\theoremstyle{remark}

\begin{document}

\renewcommand{\baselinestretch}{1.2}

\thispagestyle{empty}
\title[Pointed Hopf superalgebras]{On pointed Hopf superalgebras}
\author[andruskiewitsch, angiono and yamane]{Nicol\'{a}s Andruskiewitsch}
\author[]{Iv\'an Angiono}
\address{Facultad de Matem\'{a}tica, Astronom\'{i}a y f\'{i}sica\\
Universidad Nacional de C\'{o}rdoba \\ CIEM - CONICET,
(5000) Ciudad Universitaria \\
C\'{o}rdoba \\Argentina}
\email{andrus@mate.uncor.edu}\email{angiono@mate.uncor.edu}
\author[]{Hiroyuki Yamane}
\address{Department of Pure and Applied Mathematics,
Graduate School of Information Science and Technology,
Osaka University,
Toyonaka, Osaka, 560-0043, Japan}
\email{yamane@ist.osaka-u.ac.jp}

\thanks{Part of this work was done during the visit of the third author to the University of C\'ordoba in September 2008, supported
through Japan's Grand-in-Aid for Scientific Research (C), 19540027. The first and second authors were partially supported by ANPCyT-Foncyt, CONICET, Ministerio de Ciencia y Tecnología (Córdoba) and Secyt-UNC.}

\begin{abstract}
We discuss the relationship between Hopf superalgebras and Hopf algebras. We list the braided vector spaces of diagonal type with generalized root system
of super type and give the defining relations of the corresponding Nichols algebras.
\end{abstract}
\dedicatory{Dedicado a Jorge Vargas en su sexag\'esimo cumplea\~nos.}

\maketitle

\section*{Introduction}

The motivation for this paper is the following: in Heckenberger's classification of Nichols algebras of diagonal type with finite root system \cite{He},
there is a large class of examples that would correspond to contragredient Lie superalgebras. We want to understand this correspondence.
In fact, the explanation is very simple: Let $A$ be a Hopf algebra with bijective antipode. There is a functor, discovered by
Radford \cite{Ra} and interpreted in categorical terms by Majid \cite{Mj}, from the category of Hopf algebras in the braided category of Yetter-Drinfeld
modules over $A$ to the category of Hopf algebras with split projection to $A$: $$R \rightsquigarrow R\# A.$$ See Subsection  \ref{subsect:super-boson}
for details; $R\# A$ is called the bosonization of $R$.
If $A$ is the group algebra of $\Z/2$, then the category of super vector spaces fully embeds into the category of
Yetter-Drinfeld modules over $A$; thus, there is a functor from the category of Hopf superalgebras to the category of Hopf algebras that we denote
$H \rightsquigarrow H^{\sigma} = H\# \ku \Z /2$. This functor explains why quantum supergroups appear in the theory of Hopf algebras, as already noticed by Majid, see \cite[Chapter 10.1]{Mj2}. Majid emphasized that the theory of Lie superalgebras and Hopf superalgebras can be
reduced to the classical case using the bosonization by $\ku \Z/2$. An exposition of these facts is
given in Section \ref{sec:super}; we assume that the reader has some familiarity with the Lifting Method for the classification of Hopf algebras, see
\cite{AS-cambr}. However, the main features of this method can be read from the exposition below, dropping the signs everywhere.

Section \ref{sec:super-nichols} is devoted to Nichols algebras of diagonal type that would correspond to contragredient Lie superalgebras. We list the
related diagonal braidings and discuss the presentation by generators and relations of the corresponding Nichols algebras applying the method recently
presented  in \cite{Ang-rel}. In this way, we recover results from \cite{Y} and give a partial answer to \cite[Question 5.9]{A}. We stress that Sections 1 and 2 are independent of each other.

\subsection*{Notation}\label{subsec:notation}

Let $\ku$ be a field with char $\ku \neq 2$; all vector spaces, algebras, tensor products, etc. are over $\ku$ except
when explicitly stated. For each $N > 0$, $\G_N$ denotes the group of $N$-th roots
of 1 in $\ku$.

\medbreak  We use Sweedler's notation for the comultiplication of
a coalgebra $D$: If $x\in D$, then $\Delta(x) = x\_{1}\otimes
x\_{2}$; if $V$ is a left $D$-comodule with coaction $\lambda$ and
$v\in V$, then $\lambda(v) = v\_{-1}\ot v\_0$. If $K$ is a Hopf
algebra, then $\ydk$ denotes the category of Yetter-Drinfeld
modules over $K$.

\section{Hopf superalgebras}\label{sec:super}

\subsection{Super vector spaces}\label{subsec:super-vs}

A super vector space is a vector space $V$ graded by $\Z/2$. We
shall write $V = V^0 \oplus V^1$, to avoid confusion with the
coradical filtration. If $v\in V^j$, then we say that $v$ is
homogeneous and write $j = \vert v\vert$. A super linear map is
just a linear map between super vector spaces preserving the
$\Z/2$-grading. The tensor product of two super vector spaces $V$ and $W$
is again a super vector space, with $(V\otimes W)^i = \oplus_j V^j\otimes W^{i+j}$, $i =0, 1$.
The category of super vector spaces is symmetric, with
super symmetry $\supertau: V\ot W \to W\ot V$, $\supertau(v\ot w) = (-1)^{\vert
v\vert\vert w\vert} w\ot v$ when $v$ and $w$ are homogeneous. A
super vector space $V = V^0 \oplus V^1$ is naturally a
Yetter-Drinfeld module over the group algebra $K = \ku\Z/2$;
namely, if $\sigma$ denotes the generator of $\Z/2$, then we
define $\sigma\cdot v = (-1)^{|v|}v$ and $\delta(v) =
\sigma^{|v|}\ot v$ for $v$ homogeneous. The natural embedding of
the category of super vector spaces into $\ydk$ preserves the
braiding.

\begin{obs}\label{obs:yd-vs-super} The category $\ydk$ is semisimple and its irreducible objects
are of the form $\ku_g^{\chi}$, $g\in \Z/2$, $\chi\in \widehat{\Z/2} = \{\varepsilon, \sgn \}$, meaning that $g$ defines the coaction and $\chi$ the action on the one-dimensional vector space $\ku_g^{\chi}$. Then the category $\sv$ of super vector spaces can be identified with the subcategory of $\ydk$ whose objects have non-trivial isotypical components only of the types $\ku_e^{\varepsilon}$ (even part) or $\ku_\sigma^{\sgn}$ (odd part).
Let $\T$ be the analogous subcategory whose objects have non-trivial isotypical components only of the types $\ku_e^{\sgn}$
or $\ku_\sigma^{\varepsilon}$. Then  $\sv\ot \T \hookrightarrow \T$ (that is, $\T$ becomes a module category over $\sv$)
and $\T\ot \T \hookrightarrow \sv$. In other words, $\ydk$ is a super tensor category, with even part $\sv$ and odd part $\T$.
\end{obs}

\subsection{Superalgebras}\label{subsec:super-alg}
A \emph{superalgebra} is a $\Z/2$-graded algebra, namely
an associative algebra $A$ with a $\Z/2$-grading $A = A^0\oplus
A^1$ such that $A^iA^j\subseteq A^{i+j}$, for $0\le i,j\le 1$.
Given an associative algebra $A$, a super structure on $A$ is equivalent to an algebra automorphism of order 2, that we call $\sigma$ by abuse of notation.
Thus $\sigma$ preserves the Jacobson radical and its powers.
Let $V = V^0 \oplus V^1$ be a super vector space. Then $\End V$ is
a superalgebra, with respect to the grading $\End V = (\End V)^0
\oplus (\End V)^1$, where $(\End V)^i = \{T\in \End V:
T(V^j)\subseteq V^{i+j}\}$.

\medbreak Let $A = A^0 \oplus A^1$ be a superalgebra. Then a super
representation on a super vector space $V$ is a morphism of
superalgebras $\rho: A \to \End V$. This amounts to the same as a
module action $A\ot V \to V$ that respects the grading-- in words,
$V$ is a (left) \emph{supermodule}. A \emph{superbimodule} over $A$ is a
bimodule such that both the left and the right actions are super
linear. The corresponding categories are denoted $\supa$,
$\supaa$; morphisms are super linear and preserve both
structures.

\begin{obs}\label{obs:supermodules} Let $A^{\sigma} = A \# \ku\langle \sigma \rangle$, the smash product algebra;
as a vector space, this is $A\ot \ku\Z/2 \simeq A \oplus A\sigma$, and the
multiplication is given by
\begin{equation}\label{eqn:smash-product}
(a \# \sigma^k)(b \# \sigma^l) = (-1)^{k|b|} ab \#  \sigma^{k+l}, \qquad a,b\in A, k,l \in \{0,1\}.
\end{equation}
Then the category $\supa$ is naturally equivalent to the category $\mas$
of modules over $A^{\sigma}$. Explicitly, given an $A$-supermodule $V$ we define the action of $A^{\sigma}$ as
$$ (a \# \sigma^k) \cdot v:= (-1)^{|v|k} a \cdot v, \qquad v  \in V, a \in A, \quad k \in \Z/2. $$
Reciprocally, any object of $\mas$, with the grading defined by the action of $\sigma$ and the action of $A$ given by restriction, is an $A$-supermodule.
\end{obs}

\subsection{Supercoalgebras}\label{subsec:super-coalg}
A \emph{supercoalgebra} is a $\Z/2$-graded coalgebra,
namely an associative coalgebra $C$ (with comultiplication
$\Delta$ and counit $\varepsilon$), provided with a $\Z/2$-grading
$C = C^0\oplus C^1$ such that $\Delta(C^i)\subseteq \sum_j C^j\ot
C^{i-j}$, for $0\le i\le 1$. We write the comultiplication of $C$
by the following variation of Sweedler's notation: If $c\in C$,
then $\Delta(c) = c\^{1}\otimes c\^{2}$, where $c\^{1}$ and
$c\^{2}$ are homogeneous.

The dual of a finite-dimensional superalgebra is a supercoalgebra; thus, the matrix coalgebra $\End^c V$
is a supercoalgebra, if $V$ is a finite-dimensional super vector space.

Given an associative coalgebra $C$, a super structure on $C$ is equivalent to a coalgebra automorphism of order 2.
In particular, the coradical $C_0$, hence all terms
$C_n$ of the coalgebra filtration, are stable under $\sigma$,
\emph{i.~e.}, are super vector subspaces of $C$. Therefore the
associated graded coalgebra $\gr C = \oplus_{n\ge 0} C_n/C_{n-1}$
is a  $\Z$-graded supercoalgebra (where $C_{-1} = 0$).

An element $c\in C$ is a group-like if $c\neq 0$ and $\Delta(c) = c\ot c$; the set of all group-likes is denoted $G(C)$.
Since $\sigma$ acts by coalgebra automorphisms, it preserves $G(C)$.
We shall say that $C$ is \emph{pointed} if the coradical $C_0$ equals the
linear span of $G(C)$.
In the particular case when $C$ is the linear span of $G(C)$, a structure of super coalgebra is determined by an involution of the set $G(C)$.

\medbreak Let $C = C^0 \oplus C^1$ be a supercoalgebra. Then a (left)
\emph{supercomodule} is a super vector space $V$ provided with a
coaction $\delta:V \to C \ot V $ that is super; similarly for right supercomodules and
superbicomodules. The corresponding categories are $\supc$, $\supcr$, $\supc^C$.

\begin{obs}\label{obs:supercomodules} Let $C^{\sigma} = C \# \ku\langle \sigma \rangle$, the smash product coalgebra;
as a vector space, this is $C\ot \ku\Z/2$, and the comultiplication is given by
\begin{equation}\label{eqn:smash-coproduct} \Delta(c\# \sigma^k) = c^{(1)} \# \sigma^{\vert c^{(2)}\vert + k}  \otimes c^{(2)}\#\sigma^k,
\qquad c \in C, \quad k \in \Z/2.
\end{equation}
Then $(\gr C)^{\sigma} \simeq \gr (C^{\sigma})$.

The category $\supc$ is naturally equivalent to the category $\cocs$ of left comodules over $C^{\sigma}$.
In this case, given $V\in \supc$, the coaction $\delta: V \rightarrow C^{\sigma} \otimes V$ is defined by
$$ \delta(v) :=  v^{(-1)} \# \sigma^{\vert v ^{(0)} \vert} \otimes v^{(0)}, \qquad v\in V. $$
Reciprocally, any object of $\cocs$, with the grading  and the action of $A$ given by corestriction, is an $A$-supercomodule.
\end{obs}

\subsection{Hopf Superalgebras}\label{subsec:super-hopf}
A  \emph{superbialgebra} is a bialgebra in the category of super vector spaces, that is a $\Z/2$-graded algebra and
coalgebra $B$ (with respect to the same grading) such that
$\Delta$ and $\varepsilon$ are multiplicative, with respect to the
product in $B\otimes B$ twisted by the
super symmetry: $(a\ot b)(c\ot d) = (-1)^{\vert b\vert\vert
c\vert} ac\ot bd$. A \emph{Hopf superalgebra} is a superbialgebra
$H$ such that the identity map has a convolution inverse $\Ss\in
\End H$; $\Ss$ is called the antipode and preserves the super
grading.

\begin{exa}\label{exa:group-superhopf}
When does a usual Hopf algebra $H$ admit a Hopf superalgebra structure? If so, $\sigma$ acts by a Hopf algebra automorphism of order 2. However, this is not enough: for, assume that $H = \ku \Gamma$ is a group algebra and let $\sigma_0$ an involution of the group $\Gamma$. If $H$ is a Hopf superalgebra with respect to the automorphism $\sigma$ of $H$ defined by $\sigma_0$, then $\sigma_0 = \id$. Indeed, let $g\in \Gamma$ and let $x_g = \frac 12 (e_g + e_{\sigma_0(g)})$,
$y_g = \frac 12 (e_g - e_{\sigma_0(g)})$. Then
\begin{align*}
\Delta (e_{g^2}) &= e_{g^2}\ot e_{g^2} = e_{g}e_{g}\ot e_{g}e_{g} = e_{g}e_{g}\ot x_{g}e_{g} + e_{g}e_{g}\ot y_{g}e_{g},
\\
\Delta (e_{g})\Delta (e_{g}) &=  (e_{g}\ot (x_g + y_g))((x_g + y_g)\ot e_{g}) = e_{g}e_{g}\ot x_{g}e_{g} + e_{g}e_{\sigma_0(g)}\ot y_{g}e_{g}.
\end{align*}
Hence $y_g =0$ and $\sigma_0(g) = g$.
The same argument shows that a group-like element in a Hopf superalgebra is even. In fact, it can be shown a semisimple Hopf superalgebra over $\kk$ such
that $\Ss^2 = \id$ is purely even \cite[Cor. 3.1.2]{AEG}.
\end{exa}

A Hopf superalgebra $H$ is, in particular, a braided Hopf
algebra in $\ydk$, and there is a Hopf algebra $\hs$, the
Radford-Majid bosonization of $H$. As a vector space, $\hs = H\ot
\ku\Z/2 \simeq H \oplus H\sigma$; the multiplication and comultiplication of $\hs$ are given by \eqref{eqn:smash-product}
and \eqref{eqn:smash-coproduct}.

\begin{obs}\label{obs:eg}
Another relation between Hopf algebras and Hopf superalgebras is given in \cite[Th. 3.1.1]{AEG}:
There is a one-to-one correspondence between
    \begin{enumerate}
        \item isomorphisms classes of pairs $(\Hc,u)$, where $\Hc$ is a Hopf algebra and $u \in \Hc$ is a group-like element such that $u^2=1$, and

\medbreak\item isomorphisms classes of pairs $(H,g)$, where $H$ is a Hopf superalgebra and $g \in H$ is a group-like element such that $g^2=1$ and $gxg^{-1}=(-1)^{|x|}x$.
    \end{enumerate}

\medbreak Explicitly, given $(\Hc,u)$, $H$ is the algebra $\Hc$ with the grading given by the adjoint action of $u$ and the  comultiplication:
$$
\Delta_{\text{super}} (h) =  \Delta_{0} (h) + (-1)^{\vert h \vert} (u\ot \id) \Delta_{1} (h), \quad h\in H.
$$
Here $\Delta(h) =  \Delta_{0} (h) + \Delta_{1} (h)$, with $\Delta_{k} (h) \in H \ot H^k$, $k\in \Z/2$. Also, $g = u$.
Conversely, $\Hc = H^{\sigma} /(\sigma u - 1)$. Eventually, this correspondence leads to the classification of all finite-dimensio\-nal triangular Hopf algebras over an algebraically closed field of characteristic 0 \cite{eg}.
\end{obs}

\medbreak  Let $H = H^0 \oplus H^1$ be a Hopf superalgebra. Then
the category of supermodules $\supre$ is a tensor one (with the underlying tensor product of super vector spaces); if $V, W\in
\supre$, $v\in V$ is homogeneous, $w\in W$ and $h\in H$, then
\begin{equation}\label{eqn:supertensor}
h\cdot (v\ot w) = (-1)^{\vert h\^2\vert \vert v\vert}h\^1\cdot v
\ot h\^2\cdot w.
\end{equation}
Moreover, the equivalence $\supa \simeq \mas$ in Remark \ref{obs:supermodules} is monoidal.

Analogously, the category $\supcomh$ of supercomodules over a Hopf superalgebra $H$
is a tensor category; here, if $V, W\in
\supcomh$, $v\in V$ is homogeneous, $w\in W$, then
\begin{equation}\label{eqn:supertensor-com}
\delta (v\ot w) = (-1)^{\vert v\^0\vert \vert w\^{-1}\vert} v\^{-1} w\^{-1}
\ot v\^0 \ot w\^{0}.
\end{equation}
The equivalence $\supc\simeq \cocs$ in Remark \ref{obs:supercomodules} is also monoidal.

\medbreak  A \emph{quasitriangular Hopf superalgebra} is a pair $(H, \R)$
where $H$ is a Hopf superalgebra and $\R\in H\ot H$ is even,
invertible, and satisfies the same axioms as in the non-super case, see \emph{e.~g.} \cite{AEG}.
In particular, it provides $\supre$ with a braiding: if $V, W\in
\supre$, $v\in V$ and $w\in W$, then $c_{V,W}:
V\ot W \to W\ot V$ is given by
\begin{equation}\label{eqn:superbraiding}
c_{V,W} (v\ot w) = (-1)^{\vert v\vert \vert w\vert}\R\cdot (w \ot
v).
\end{equation}
The element $\R$ is called a universal $R$-matrix.
Furthermore,  any braiding in $\supre$ arises from
a universal $R$-matrix, cf. \cite[10.4.2]{Mon}.
Also, $(H, \R)$ is a \emph{triangular Hopf superalgebra} if $\supre$ is symmetric for the previous braiding.

\subsection{Hopf supermodules, Hopf superbimodules and Yetter-Drinfeld supermodules}\label{subsect:super-yd}

Let $H$ be a Hopf superalgebra with bijective antipode. There is a
hierarchy of special modules over $H$; the proofs of the statements below are
adaptations of the usual proofs for Hopf algebras. We leave to the reader the pleasant task of
checking that the signs match. Analogous results in the more general context of braided categories have been proved in \cite{Bes, BD}.

\medbreak\noindent $\bullet$ A \emph{Hopf supermodule} over $H$ is a super vector space $V$ that
    is simultaneously a supermodule and a supercomodule, with
    compatibility saying that the coaction $\delta:V \to H\ot V$
    is morphism of $H$-supermodules. If $U$ is a super vector space,
then $H\ot U$ is a Hopf supermodule over $H$ with the  action and coaction on the left.
If $V$ is a supercomodule, then set $V^{\co H}= \{v\in V: \delta (v) = 1\ot v\}$.
There is a \emph{Fundamental theorem for Hopf supermodules:} the category of Hopf supermodules over $H$ is equivalent to the category of super vector spaces,
via $V\mapsto V^{\co H}$. Explicitly, if $V$ is a Hopf supermodule, then the multiplication $\mu:H\ot V^{\co H} \to V$ is an
    isomorphism of Hopf supermodules.

\medbreak
    \noindent $\bullet$ A \emph{Hopf superbimodule} is a super vector space $V$ that
    is simultaneously a superbimodule and a superbicomodule, with
    compatibility saying that both coactions $\lambda:V \to H\ot
    V$ and $\rho:V \to V\ot H$ are morphisms of $H$-superbimodules. The category $\suphbim$ of Hopf superbimodules is a tensor one, with tensor product $\ot_H$.

\medbreak
    \noindent $\bullet$ A \emph{Yetter-Drinfeld supermodule} over $H$ is a super vector space $V$ that
    is simultaneously a supermodule and a supercomodule, with
    compatibility saying that
\begin{equation}\label{eqn:superyd}
\delta(h\cdot v) = (-1)^{|v^{(-1)}| (|h^{(2)}|+|h^{(3)}|)+|h^{(2)}||h^{(3)}|}h^{(1)}v^{(-1)}\Ss(h^{(3)}) \otimes h^{(2)}\cdot v^{(0)}.
\end{equation}

The category $\ydsh$ of Yetter-Drinfeld supermodules is tensor equivalent to $\suphbim$. Explicitly,
\begin{align*}
M\in \suphbim &\rightsquigarrow  V = M^{\co H} = \{m\in M: \rho(m) = m \ot 1\},
\end{align*} with action and coaction $h\cdot v = (-1)^{|v||h^{(2)}|} h^{(1)}v\Ss(h^{(2)})$, $\delta = \lambda$;
\begin{align*}
V\in \ydsh &\rightsquigarrow M = V \ot H, &
x\cdot (v\ot h) y &= (-1)^{|v||x^{(2)}|} x^{(1)}v\ot x^{(2)}hy,
\\ \lambda (v\ot h) &= v^{(-1)}h^{(1)} \otimes ( v^{(0)}\ot h^{(2)}), & \rho (v\ot h) &= (v \otimes h^{(1)}) \ot h^{(2)}
\end{align*}
for $v\in V$, $h,x,y\in H$.

\medbreak
    \noindent $\bullet$ The tensor category $\suphbim$ of Hopf superbimodules is braided, with braiding
$c_{M,N}: M\otimes_H N \rightarrow N \otimes_H M$, $M, N\in \suphbim$, given by
    $$ c_{M,N}(m\otimes n) = (-1)^{(|m^{(0)}| + |m^{(-1)}|)(|n^{(0)}| + |n^{(1)}|)} m^{-2}  n ^{(0)} \Ss(n ^{(1)})\Ss(m ^{(-1)}) \otimes m^{(0)}n ^{(2)}. $$
Thus, $\ydsh$ is braided, with braiding $c_{V,W}: V\otimes W \rightarrow W \otimes V$, $V, W\in \ydsh$, given by
\begin{equation}\label{eqn:super-trenza}
c_{X,Y}(x\otimes y) = (-1)^{|x^{(0)}||y|} x^{(-1)} \cdot y \otimes x^{(0)}.
\end{equation}

\begin{obs}\label{obs:YetterDrinfeldsupermodules}
For each Hopf superalgebra $H$ there exists a full embedding of braided tensor categories $\incl: \ydsh \hookrightarrow
{}^{H^{\sigma}} _{H^{\sigma}} \mathcal{YD}$, given by the restriction of equivalences in Remarks \ref{obs:supermodules} and \ref{obs:supercomodules}.
\end{obs}

\begin{exa}\label{exa:hopf-even}
Let $H$ be a purely even Hopf superalgebra, that is a usual Hopf algebra with trivial grading. If $V\in \ydh$
and $k\in \Z/2$, then $V[k] = V$ with all elements of degree $k$, is an object in $\ydsh$. Moreover, if
$V$ is irreducible in $\ydh$, then $V[k]$ is irreducible in  $\ydsh$. We claim that any irreducible in  $\ydsh$ is of the form $V[k]$ as above.
\end{exa}

\pf
If $W$ is an irreducible module in $\ydsh$, then $W^0$ and $W^1$ are sub-objects in $\ydsh$, hence $W = W^k$ for some $k$.
Thus $W = U[k]$, for $U$ an irreducible sub-object of $W^k$ in $\ydh$. \epf

\subsection{Hopf superalgebras with projection and bosonization}\label{subsect:super-boson}

In this subsection we consider bosonization and Hopf superalgebras with projections; we note that this construction can be done for general
braided categories, see \cite{Bes,BD}.

\medbreak Let $H$ be a Hopf superalgebra with bijective antipode. If $R$ is a Hopf algebra in the braided category $\ydsh$, then we have a Hopf superalgebra $R\#H$: it has $R \otimes H$ as underlying super vector space, and its structure is defined by
\begin{align*}
    (a\#h)(b\#f) & := (-1)^{|h^{(2)}||b|} a(h^{(1)} \cdot b) \# h^{(2)} f, \qquad
    1  := 1_R \# 1_H, \\
    \Delta(a\#h) & := (-1)^{|(a^{(2)})_{(0)}||h^{(1)}|} a^{(1)}\# (a^{(2)})_{(-1)}(h^{(1)} \otimes (a^{(2)})_{(0)} \# h^{(2)}, \\
    \varepsilon(a\#h) & := \varepsilon_R(a) \varepsilon_H( h), \\
    \Ss(a\#h) & := (-1)^{|a_{(0)}||h|} \left( 1\# \Ss_H(a_{(-1)}h) \right) \left(\Ss_R(a_{(0)}) \# 1 \right).
\end{align*}
for each $a,b \in R$ and $h,f \in H$.
By Remark \ref{obs:YetterDrinfeldsupermodules}, the image of $R$ under the full embedding $\incl$
is a Hopf algebra in $^{H^{\sigma}} _{H^{\sigma}} \mathcal{YD}$. It is straightforward to prove that
\begin{equation}
(R\# H)^{\sigma} \cong \incl (R) \# H ^{\sigma}.
\end{equation}

Let $\iota: H \hookrightarrow L$ and $\pi: L \twoheadrightarrow H$ be morphisms of Hopf superalgebras  satisfying $\pi \circ \iota = \id_H$.
Consider the subalgebra of coinvariants
$$ R:= L^{\coin H} = \{x \in L: \, (\id \otimes \pi) \Delta(x)= x \otimes 1 \}. $$
This is a Hopf algebra in the category $^{H} _{H} \mathcal{YDS}$ and there exists an isomorphism of Hopf superalgebras $L \cong R \# H$.
Now, $\iota$ and $\pi$ induce Hopf algebra morphisms $\iota_{\sigma}: H^{\sigma} \hookrightarrow L^{\sigma}$ and $\pi_{\sigma}: L^{\sigma} \twoheadrightarrow H^{\sigma}$ such that $\pi_{\sigma} \circ \iota_{\sigma} = \id_{H^{\sigma}}$. Then $\incl(R)$ coincides with the subalgebra of coinvariants
$(L^{\sigma})^{\coin H^{\sigma}}$, see \emph{e.~g.} \cite[Lemma 3.1]{AHS}.


\subsection{Nichols superalgebras}\label{subsec:super-nichols}
Let $H$ be a Hopf superalgebra with bijective antipode. The constructions and results of \cite[Section 2]{Sbg} hold in
the braided abelian category $^{H} _{H} \mathcal{YDS}$. We summarize:

\begin{prop} \cite[Section 2]{Sbg}. Let $V$ be a Yetter-Drinfeld supermodule over $H$. Then there is a unique (up to isomorphisms) graded
Hopf algebra $\toba(V) = \oplus_{n\in\N_0}\toba^n(V)$ in $\ydsh$ with the following properties:

\begin{itemize}
  \item $\toba^0(V) \simeq\ku$,
  \item $V \simeq \toba^1(V) = \mathcal P(\toba(V))$ (the space of primitive elements),
  \item $\toba^1(V)$ generates the algebra $\toba(V)$.
\end{itemize}
\end{prop}

Explicitly, $\toba(V) \simeq T(V) /\J(V)$, where the ideal $\J(V) = \oplus_{n\ge 2} \J^n(V)$ has homogeneous components $\J^n(V)$
that equal the kernel of the quantum symmetrizer $\ker \mathfrak S_{n}$
\cite{Wo, Sbg}. To be more precise, the braid group in $n$ letters $\mathbb B_{n}$ acts on $T^n(V)$ via the braiding in $\ydsh$.
Let $\pi: \mathbb B_{n} \to \mathbb S_{n}$ be the
natural projection and let $s:\mathbb S_{n} \to \mathbb B_{n}$ be the so-called Matsumoto (set-theoretical) section. The element
$\mathfrak S_{n} := \sum_{\sigma \in \mathbb S_{n}} s(\sigma)$ of the group algebra of $\mathbb B_n$ is called the quantum symmetrizer;
it acts on $T^n(V)$ and its kernel is $\J^n(V)$.
Because of this explicit description, we conclude from Remark \ref{obs:YetterDrinfeldsupermodules} that the Nichols algebra functor commutes with the full embedding $\incl$:
\begin{equation}\label{eqn:embedding-nichols}
\toba\left(\incl (V)\right) \simeq \incl\left(\toba(V)\right).
\end{equation}

\subsection{The lifting method for Hopf superalgebras}\label{subsec:super-lifting}
Let $H$ be a Hopf superalgebra with bijective antipode. If the coradical $H_0$ of $H$ is a Hopf sub-superalgebra, then the coradical filtration is also an algebra filtration and the associated graded coalgebra $\gr H$ is a graded Hopf superalgebra. Furthermore,
the homogeneous projection $\gr H \to H_0$ splits the inclusion, hence gives rise to a decomposition $\gr H \simeq R \# H_0$. The graded Hopf algebra $R  = \oplus_{n\ge 0} R^n\in \ydsho$ has the following properties:

\begin{itemize}
  \item $R^0 \simeq \ku$,
  \item $R^1 = \mathcal P(R)$.
\end{itemize}
Thus the subalgebra generated by $V:= R^1$ is isomorphic to the Nichols algebra $\toba(V)$. In this way, the Lifting Method \cite{AS-cambr} can be adapted to the setting of Hopf superalgebras whose coradical is a Hopf sub-superalgebra.
However, there is no need to start over again since classification problems of Hopf superalgebras reduce to analogous
classification problems of Hopf algebras via the functor $H \rightsquigarrow H^{\sigma}$.
This principle is illustrated by the following facts:

\medbreak\noindent $\bullet$ The coradical of $\hs$ is $H_0 \oplus H_0\sigma$ and $G(H^{\sigma}) = G(H) \times \langle\sigma\rangle$.
More generally, the coradical filtration of $\hs$ is $\hs_n = H_n \oplus H_n\sigma$.

\medbreak\noindent $\bullet$ $\hs$ is pointed if and only if $H$ is pointed.

\medbreak\noindent $\bullet$ The coradical of $\hs$ is a Hopf subalgebra if and only if
the coradical of $H$ is a Hopf sub-superalgebra. If this is the case, then  $$
\gr (H^{\sigma}) \simeq (\gr H)^{\sigma} \simeq (R \# H_0)^{\sigma} \simeq \incl (R) \# (H_0)^{\sigma}.$$

\medbreak\noindent $\bullet$ The Hopf algebra  $H$ is generated (as algebra)
by group-like and skew-primitive elements (generated in degree one, for short) if and only if
$\hs$ is generated in degree one.

\begin{obs}\label{obs:gen-deg-one} It was conjectured that a finite-dimensional pointed Hopf algebra over
$\ku$ is generated in degree one \cite[1.4]{AS}.
This Conjecture was verified in various cases, see \emph{e.~g.} \cite[7.6]{AS-ann-ec-sci}, \cite{GG}, \cite[2.7]{Ang-GI}, \cite[4.3]{Ang-rel}.
The validity of the conjecture would imply the validity of the analogous one for Hopf superalgebras.
\end{obs}

Because of these considerations, we see that the theory of Hopf superalgebras is naturally a part of the theory of Hopf algebras.

\begin{exa}\label{exa:clasification-Z2}
There is a full embedding from the category of Lie superalgebras to the category of pointed Hopf algebras with group $\Z/2$,
given by $\mathfrak g \rightsquigarrow U(\mathfrak g)^{\sigma}$.
\end{exa}

\medbreak In particular, we see that the classification of \fd{} pointed Hopf superalgebras $H$ with a fixed group of group-like elements $\Gamma$ reduces to
the classification of \fd{} pointed Hopf algebras $K$ such that
\begin{itemize}
  \item[$\circ$] $G(K) \simeq \Gamma \times \Z/2$,
  \item[$\circ$] there exists a projection of Hopf algebras $K \to \ku \Z/2$ that splits the inclusion (from the second factor above).
\end{itemize}

\begin{exa}\label{exa:superyd-gamma-abelian}
Let $\Gamma$ be a finite abelian group. Assume that $\ku$ is algebraically closed. Then any irreducible object in $\ydg$ has dimension one and is of the form $\ku_g^{\chi}$, $g\in \Gamma$, $\chi\in \widehat{\Gamma}$, where $g$ determines the coaction and $\chi$ the action.  By Example \ref{exa:hopf-even}, any irreducible object in $\ydg$
is of the form $\ku_g^{\chi}[k]$, $g\in \Gamma$, $\chi\in \widehat{\Gamma}$, $k\in \Z/2$.
The corresponding isotypical component of $V\in \ydg$ is denoted $V_g^{\chi}[k]$. Thus, any \fd{} $V\in \ydsg$ has a basis $x_1, \dots, x_{\theta}$ with
$x_j \in V_{g_j}^{\chi_j}[k_j]$, $g_j\in \Gamma$, $\chi_j\in \widehat{\Gamma}$, $k_j = \vert x_j\vert\in \Z/2$.
\end{exa}

\begin{prop} \label{prop:trenza-diagonal} Let $g_j\in \Gamma$, $\chi_j\in \widehat{\Gamma}$, $k_j\in \Z/2$, $1\le j\le \theta$.
For $1\le i,j\le \theta$, set $q_{ij} := \chi_j(g_i)$ and
\begin{equation}\label{eqn:super-diagonal}
\widetilde q_{ij} = \begin{cases} q_{ij}, & i\neq j,
\\ (-1)^{k_i} q_{ii}, & i=j. \end{cases}
\end{equation}
Let $V\in \ydsg$ with a basis $x_1, \dots, x_{\theta}$, such that $x_j \in V_{g_j}^{\chi_j}[k_j]$.
Then the Nichols superalgebra $\toba(V)$ has finite dimension if and only if the connected components of the generalized Dynkin diagram corresponding to the matrix $(\widetilde q_{ij})_{1\le i,j\le \theta}$ belong to the list in \cite{He}. \end{prop}

\pf By \eqref{eqn:super-trenza} and \eqref{eqn:embedding-nichols}, we are reduced to consider the Nichols algebra of the braided vector space of diagonal type with matrix $\left((-1)^{\vert x_i\vert\vert x_j\vert} q_{ij}\right)_{1\le i,j\le \theta}$. Now this matrix and $(\widetilde q_{ij})_{1\le i,j\le \theta}$ are twist-equivalent \cite[Def. 3.8]{AS-cambr}, hence their Nichols algebras have the same dimension \cite[Prop. 3.9]{AS-cambr}.
\epf

\section{Generalized root systems and Nichols algebras}\label{sec:super-nichols}

\subsection{Generalized root systems}\label{subsec:gen-root-syst}

We recall now the generalization of the notion of a root system given in \cite{HY}.

Fix two non-empty sets $\cX$ and $I$, where $I$ is finite, and denote by $\{\alpha_i \}_{i \in I}$ the canonical basis of $\Z^I$.

\begin{definition}{\cite{HY,CH1}}
Assume that for each $i \in I$ there exists a map $r_i: \cX \rightarrow \cX$, and for each $X \in \cX$ a generalized Cartan matrix
$A^X= (a^X_{ij})_{i,j \in I}$ in the sense of \cite{K} satisfying
\begin{enumerate}
  \item for all $i \in I$, $r_i^2=id$, and
  \item for all $X \in \cX$ and $i,j \in I$: $a^X_{ij}=a^{r_i(X)}_{ij}$.
\end{enumerate}
We say that the quadruple $\cC:= \cC(I, \cX, (r_i)_{i \in I}, (A^X)_{X \in \cC})$ is a \emph{Cartan scheme}.

Given $i \in I$ and $X \in \cX$, $s_i^X$ denotes the automorphism of $\Z^I$ such that
$$ s_i^X(\alpha_j)=\alpha_j-a_{ij}^X\alpha_i, \qquad j \in I. $$

The \emph{Weyl groupoid} of $\cC$ is the groupoid $\cW(\cC)$ for which:
\begin{enumerate}
 \item the set of objects is $\cX$, and
 \item the morphisms are generated by $s_i^X$, if we consider $s_i^X \in \Hom(X, r_i(X))$, $i \in I$, $X \in \cX$.
\end{enumerate}
Each morphism $w \in \Hom(\cW,X_1)$ is a composition $s_{i_1}^{X_1}s_{i_2}^{X_2} \cdots s_{i_m}^{X_m}$,
where $X_j=r_{i_{j-1}} \cdots r_{i_1}(X_1)$, $i \geq 2$. We shall write $w= \id_{X_1} s_{i_1} \cdots s_{i_m}$ to indicate that $w \in \Hom(\cW,X_1)$,
because the $X_j$'s are univocally determined by the first one and the sequence $i_1, \cdots, i_m$.
\end{definition}

\begin{definition}{\cite{HY,CH1}}
Given a Cartan scheme $\cC$, and for each $X \in \cX$ a set $\de^X \subset \Z^I$, define $m_{ij}^X:= |\de^X \cap (\N_0\alpha_i+\N_0 \alpha_j)|$. We say that
 $\cR:= \cR(\cC, (\de^X)_{X \in \cX} )$ is a \emph{root system of type} $\cC$ if
\begin{enumerate}
  \item for all $X \in \cX$, $\de^X= (\de^X \cap \N_0^I) \cup  -(\de^X \cap \N_0^I)$,
  \item for all $i \in I$ and all $X \in \cX$, $s_i^X(\de^X)=\de^{r_i(X)}$,
  \item for all $i \in I$ and all $X \in \cX$, $\de^X \cap \Z \alpha_i= \{\pm \alpha_i \}$,
  \item for all $i \neq j \in I$ and all $X\in \cX$, $(r_ir_j)^{m_{ij}^X}(X)=X$.
\end{enumerate}
\end{definition}

We call $\de^X_+:= \de^X \subset \N_0^I$ the set of \emph{positive roots} of $X$, and $ \de^X_-:= - \de^X_+$ the set of \emph{negative roots}. By simplicity
we will write $\cW$ in place of $\cW(\cC)$ when $\cC$ is understood, and for any $X \in \cX$:
\begin{align}\label{defHom}
    \Hom(\cW,X) & := \cup_{Y \in \cX} \Hom(Y,X),
    \\ \label{defrealroot} \de^{X \ re} &:= \{ w(\alpha_i): \ i \in I, \ w \in \Hom(\cW,X) \}.
\end{align}
The elements of $\de^{X \ re}$ are the \emph{real roots} of $X$.

We say that $\cR$ is \emph{finite} if $\de^X$ is finite for all $X\in \cX$. In such case all the roots are real, see \cite[Prop. 2.12]{CH1}, and
for each pair $i \neq j \in I$ and each $X \in \cX$, $\alpha_i + k \alpha_j \in \de^X$ if and only if $0 \leq k \leq -a_{ij}^X$. Therefore,
\begin{equation}\label{aij}
 a_{ij}^X = -\max \{ k \in \N_0: \ \alpha_i + k \alpha_j \in \de^X \}.
\end{equation}

\begin{exa}\label{example:super root systems}
By \cite[Example 3]{HY}, the root system associated to a finite dimensional contragradient Lie superalgebra is a generalized root system in this context.
We describe them case by case, considering the irreducible root systems. We call them \emph{super root systems}.

\medbreak
\noindent \emph{Type $A_\theta$}: We need to consider a parity of the simple roots $p(\alpha_i)$, and extend it to a group homomorphism
$p: \Z^\theta \to \{\pm 1 \}$. The set $\cX$ is determined as follows: we have a symmetry $s_i$ from one point $X$ to a different one $\hat X$ if $p(\alpha_i)=-1$. The new parity function $\hat p$ is determined from $p$ and $s_i$:
$$\hat p(\alpha_k)= p( s_i(\alpha_k))= p(\alpha_k+m_{ik}\alpha_i)=p(\alpha_k)(-1)^{m_{ik}}.$$
That is, it changes the parity of the vertices $k$ which are connected to $i$, and keeps the parity of the non-connected vertices.
In this way, $\cX$ can have more than one element depending on the parity of the simple roots, but for any $X \in \cX$ we have
the same set of positive roots,
\begin{equation}\label{root system A}
\de^X_+ = \{ \ub_{ij}:= \alpha_i+\alpha_{i+1}+ \cdots + \alpha_j: \ 1 \leq i \leq j \leq \theta \}.
\end{equation}

\medbreak
\noindent \emph{Type $B_\theta$}: As above, $\cX$ can have more than a point, and again the symmetries that go from a point to a different one are symmetries
of odd vertices, with the same changes. Anyway, the set of positive roots is the same for any $X \in \cX$,
\begin{equation}\label{root system B}
\de^X_+ = \{ \ub_{ij}: \ 1 \leq i \leq j \leq \theta \} \cup \{ \vb_{ij}:= \ub_{i,\theta} + \ub_{j, \theta}: \ 1 \leq i < j \leq \theta \}.
\end{equation}

\medbreak
\noindent \emph{Types $C_\theta, D_\theta$}: As above we consider a parity function $p: \Z^\theta \to \Z_2$. Following the classical literature, there are sets $\de^{X(C)}_+$ of type $C$ and sets $\de^{X(D)}_+$ of type $D$,
described as follows:
\begin{align} \label{root system C}
\de^{X(C)}_+ = & \ \{ \ub_{ij}: \ 1 \leq i \leq j \leq \theta \}
\\ & \cup \{ \wb_{ij}:= \ub_{i,\theta} + \ub_{j, \theta-1}: \ 1 \leq i < j \leq \theta-1 \}  \nonumber
\\ & \cup \{ \widetilde \wb_{i}:= \ub_{i,\theta-1} + \ub_{i, \theta}: \ 1 \leq i \leq \theta-1, \ p(\ub_{i, \theta-1})=1 \}, \nonumber
\end{align}
\begin{align} \label{root system D}
\de^{X(D)}_+ = & \ \{ \ub_{ij}: \ 1 \leq i \leq j \leq \theta, (i,j) \neq ( \theta-1, \theta) \}
\\ & \cup \{ \alpha_{\theta-1}+ \alpha_ \theta: p(\alpha_{\theta-1})=-1 \} \nonumber
\\ & \cup \{ \widetilde \ub_{i}:= \ub_{i,\theta-2} + \alpha_{\theta}: \ 1 \leq i \leq \theta-2 \} \nonumber
\\ & \cup \{ \zb_{ij}:= \ub_{i,\theta} + \ub_{j, \theta-2}: \ 1 \leq i < j \leq \theta-2 \}   \nonumber
\\ & \cup \{ \widetilde \zb_{i}:= \ub_{i,\theta} + \ub_{i, \theta-2}: \ 1 \leq i \leq \theta-2, \ p(\ub_{i, \theta-1})=-1 \}. \nonumber
\end{align}

\medbreak
\noindent \emph{Type $D(2,1;\alpha)$}: We have four possible sets of roots,
\begin{align}\label{root system D(2,1,alpha)1}
\de^{X_0}_+  = & \ \{ \alpha_1, \alpha_2, \alpha_3, \alpha_1+\alpha_2, \alpha_1+ \alpha_3, \alpha_2+ \alpha_3, \alpha_1+ \alpha_2+ \alpha_3 \},
\\ \label{root system D(2,1,alpha)2}
\de^{X_k}_+ = & \ \{ \alpha_1, \alpha_2, \alpha_3, \alpha_1+ \alpha_2+ \alpha_3, \alpha_1+ \alpha_2+ \alpha_3 + \alpha_k \}
\\ & \cup \{ \alpha_k+ \alpha_j: \ j \in \{1,2,3 \} \setminus \{ k \} \}, \nonumber
\end{align}
where $k \in \{1,2,3\}$. Here, $s_k(\de^{X_0})= \de^{X_k}$.

\medbreak
\noindent \emph{Type $F(4)$}: in this case $|\cX|=6$. One of the sets of roots is
\begin{align}\label{root system F4}
\de^{X}_+  = & \ \{ \alpha_1, \ \alpha_1+\alpha_2, \ \alpha_1+ \alpha_2+ \alpha_3, \ \alpha_1+ \alpha_2+ 2\alpha_3, \ \alpha_1+ 2 \alpha_2+ \alpha_3,
\\ &  \alpha_1+ \alpha_2+ \alpha_3+\alpha_4, \ \alpha_1+ \alpha_2+ 2 \alpha_3+\alpha_4, \ \alpha_1+ 2 \alpha_2+ 2 \alpha_3+\alpha_4, \nonumber
\\ & \alpha_1+2 \alpha_2+3 \alpha_3+ 2\alpha_4, \ \alpha_2+ \alpha_3+\alpha_4, \  \alpha_2+ 2 \alpha_3+\alpha_4,  \ \alpha_2+ \alpha_3, \nonumber
\\ & \alpha_2+ 2 \alpha_3, \ \alpha_2, \ \alpha_3, \ \alpha_3+ \alpha_4, \ \alpha_4 \}. \nonumber
\end{align}
The other sets of roots are obtained applying the symmetries $s_i$, once one determines $a_{ij}^X$ as in \eqref{aij}.

\medbreak
\noindent \emph{Type $G(3)$}: now, $|\cX|=4$, and one of these sets of positive roots is
\begin{align}\label{root system G3}
\de^{X}_+  =  \ \{ &\alpha_1, \ \alpha_1+\alpha_2, \ \alpha_1+ \alpha_2+ \alpha_3, \ \alpha_1+ 2\alpha_2+ \alpha_3,
\\ &  \alpha_1+ 3 \alpha_2+ \alpha_3, \ \alpha_1+ 3\alpha_2+ 2 \alpha_3, \ \alpha_1+ 4\alpha_2+ 2 \alpha_3, \nonumber
\\ &  \alpha_2,  \ \alpha_2+  \alpha_3, \ 2 \alpha_2+  \alpha_3, \ 3\alpha_2+  \alpha_3, \ 3 \alpha_2+ 2 \alpha_3,\ \alpha_3\}. \nonumber
\end{align}
We obtain the other sets of positive roots by determining $a_{ij}^X$ as in \eqref{aij} and applying the symmetries $s_i$.
\end{exa}

\bigbreak

Now we recall the definition of the Weyl groupoid attached to a braided vector space $(V,c)$ of diagonal type given in \cite{He1}, see also \cite{AA}. Fix a basis $\{ x_1, \ldots,x_{\theta}\}$ and   scalars $q_{ij} \in \ku^\times$ such that $c(x_i \ot x_j)= q_{ij} x_j \ot x_i$. Let $\chi: \zt \times \zt \to \ku^\times$ be the bilinear form such that $\chi(\alpha_i, \alpha_j)=q_{ij}$. Following \cite{He1}, $\Delta^V_+ $ denotes the set of degrees of a PBW basis of $\toba(V)$, counted with their multiplicities. It is remarked in \cite{He1} and proved in \cite{AA} that this set does not depend on the PBW basis.

For each $n \in \N$ we set the following polynomials in $q$:
$$ \binom{n}{j}_q = \frac{(n)_q!}{(k)_q! (n-k)_q!}, \quad \mbox{where }(n)_q!= \prod_{j=1}^n (k)_q, \quad \mbox{and } (k)_q= \sum_{j=0}^{k-1} q^j. $$

Let $\cX$ be the set of ordered bases of $\zt$. For each $F= \left\{ f_1, \ldots,f_{\theta} \right\} \in \cX$, set $q_{ij}^F= \chi(f_i,f_j)$. Define
\begin{equation}
a_{ij}(F):= -\min \left\{ n \in \mathbb{N}_0: (n+1)_{q^F_{ii}} (1-(q_{ii}^F)^n q^F_{ij} q^F_{ji} )=0 \right\}, \label{mij}
\end{equation}
for each $1 \leq i \neq j \leq \theta$, and set $s_{i,F} \in \Aut(\mathbb{Z}^{\theta})$ such that $s_{i,F}(f_j)=f_j-a_{ij}(F)f_i$. Here $a_{ii}=2$.

Consider for $\cG=\Aut(\mathbb{Z}^{\theta}) \times \cX$ the groupoid structure given as follows: the set of objets is $\cX$ and the morphisms are
$x \stackrel{(g,x)}{\longrightarrow} g(x)$.

Then the \emph{Weyl Groupoid} $W(\chi)$ of $\chi$ is the least subgroupoid of $\cG$ such that $(id,E) \in W(\chi)$, and if $(id,F) \in W(\chi)$ and $s_{i,F}$ is defined, then $(s_{i,F},F) \in W(\chi)$.

The generalized root system for each object $F$ is $\de^{V_F}$, where $(V_F,c_F)$ is the braided vector space of diagonal type whose braiding matrix is $(q_{ij}^F)$. It satisfies the axioms of a root system by \cite{He1}.

\begin{obs}\label{obs:simetria trenza 1}
If $a_{ij}^F=0$, then for all $k \neq i,j$,
$$q_{jj}^{s_i(F)}=q_{jj}^F, \qquad q_{jk}^{s_i(F)}q_{kj}^{s_i(F)}=q_{jk}^Fq_{kj}^F.$$
\end{obs}

\begin{obs}\label{obs:simetria trenza 2}
If $(q_{ii}^F)^{a_{ij}(F)}=q_{ij}^Fq_{ji}^F$ for all $j \neq i$, then $(q_{kj}^{s_i(F)})$ is the transpose matrix of $(q_{kj}^F)$ and then the braiding matrices are twist equivalent. In consequence, $\de^{V_{s_i(F)}}= \de^{V_F}$.
\end{obs}

\subsection{Diagonal braidings of super type}\label{subsec:braiding-super}
From now on, $\ku$ is an algebraically closed field of characteristic 0.

We shall characterize the Nichols algebras whose root system is one of those associated to a contragradient finite-dimensional Lie superalgebra.

\bigbreak
First we recall some definitions following \cite{He}. The \emph{generalized
Dynkin diagram} associated to a braided vector space of diagonal
type, with braiding matrix $(q_{ij})_{1 \leq i,j \leq \theta}$ is a
graph with $\theta$ vertices, each of them labeled with the
corresponding $q_{ii}$, and an edge between two vertices $i,j$ if $q_{ij}q_{ji} \neq 1$,
labeled with this scalar. In this way two braided vector spaces of diagonal type have the same
generalized Dynkin diagram if and only if they are twist
equivalent.

\bigbreak
A \emph{simple chain} of length $\theta$ is a braided vector space of diagonal type whose braiding matrix $(q_{ij})_{1 \leq i,j \leq \theta}$ satisfies

\begin{itemize}
\medbreak\item $(1+q_{11})(1-q_{11}q_{12}q_{21})=(1+q_{\theta \theta})(1-q_{\theta \theta}q_{\theta, \theta-1}q_{\theta-1, \theta})=0$,

\medbreak\item $q_{ij}q_{ji}=1$ if $1 <i, j<\theta$, $|i-j|>1$,

\medbreak\item for any $1 <i<\theta$, $q_{ii}=-1$, $q_{i-1,i}q_{i,i-1}q_{i+1,i}q_{i,i+1}=1$, or $q_{ii}q_{i-1,i}q_{i,i-1}=q_{ii}q_{i+1,i}q_{i,i+1}=1$.
\end{itemize}

\bigbreak
Here $C(\theta,q;i_1,\ldots ,i_j)$ denotes a simple chain such that $q=q_{\theta \theta}^2 q_{\theta, \theta-1}q_{\theta-1, \theta}$, and $q_{i-1,i}q_{i,i-1}=q$ if and only if $i \in \{i_1, \ldots, i_j \}$.

\begin{theorem}
Let $(V,c)$ a braided vector space of diagonal type, with braiding matrix $(q_{ij})$. Assume that its generalized Dynkin diagram is connected.
Then $\toba(V)$ has a super root system if and only if its generalized Dynkin diagram is one of the following ones:

\bigbreak
\noindent \textbf{Type $A_{\theta}$:}
\begin{equation}\label{tipoA}
C(\theta,q;i_1,\ldots ,i_j), \qquad \theta \in \N, q \in \ku^\times, q^2 \neq 1, 1 \leq i_1 < i_2 < \cdots < i_k \leq \theta.
\end{equation}

\bigbreak
\noindent \textbf{Type $B_{\theta}$:}
\begin{align}\label{tipoB1}
& \xymatrix{ \circ^{q} \ar@{-}[r]^{q^{-1}} & \circ^{\zeta}}, \qquad \zeta \in \G_3, \ q \in \ku \setminus \{0,1,-1,\zeta, \zeta^2 \},
\\ & \put(47,3){\oval(100,15)} C(\theta-1,q^2;i_1,\ldots ,i_j) \xymatrix{ \ar@{-}[r]^{q^{-2}} & \circ^{q} }
, \qquad \theta \in \N, \ q \in \ku^\times, q \neq \pm 1, \label{tipoB2}
\\ & \put(50,3){\oval(110,15)} C(\theta-1,- \zeta^2;i_1,\ldots ,i_j) \xymatrix{ \ar@{-}[r]^{-\zeta} & \circ^{\zeta} }
, \qquad \theta \in \N, \ \zeta \in \G_3. \label{tipoB3}
\end{align}

\bigbreak
\noindent \textbf{Type $C_{\theta}$}, $\theta \in \N$, $q \in \ku^\times$, $q^4 \neq 1$:
\begin{equation}\label{tipoC}
\put(45,3){\oval(100,15)} C(\theta-1,q;i_1,\ldots ,i_j) \xymatrix{ \ar@{-}[r]^{q^{-2}} & \circ^{q^2} }.
\end{equation}

\bigbreak
\noindent \textbf{Type $D_{\theta}$}, $\theta \in \N$, $q \in \ku^\times$, $q^2 \neq 1$:
\begin{align}\label{tipoD1}
& \put(50,3){\oval(120,15)} C(\theta-2,q^{-1};i_1,\ldots ,i_j) \xymatrix{ \ar@{-}[r]^{q} \ar@{-}[d]^{q} & \circ^{q^{-1}} \\ \circ^{q^{-1}} & },
\\ \label{tipoD2} & \put(45,3){\oval(110,15)} C(\theta-2,q;i_1,\ldots ,i_j) \xymatrix{ \ar@{-}[r]^{q^{-1}} \ar@{-}[d]^{q ^{-1}} & \circ^{-1} \ar@{-}[dl]^{q^2} \\ \circ^{-1} & }.
\end{align}

\bigbreak
\noindent \textbf{Type $D(2,1; \alpha)$}, $q,r,s \in
\kk^\times\setminus \{1\}$, $qrs=1$:
\begin{align}\label{tipoD21-1}
& \xymatrix{ \circ^{q} \ar@{-}[r]^{q^{-1}} & \circ^{-1} \ar@{-}[r]^{r^{-1}} & \circ^{r} },
\\ & \xymatrix{ & \circ^{-1} \ar@{-}[rd]^{r} & \\ \circ^{-1} \ar@{-}[ru]^{q} \ar@{-}[rr]^{s} &  & \circ^{-1} }.
\end{align}

\bigbreak
\noindent \textbf{Type $F(4)$}, $q \in \ku^\times$, $q^2, q^3 \neq 1$:

\begin{align}\label{tipoF(4)-1}
& \xymatrix{ \circ^{-1} \ar@{-}[r]^{q^{-1}} & \circ^{q} \ar@{-}[r]^{q^{-2}} & \circ^{q^2} \ar@{-}[r]^{q^{-2}} & \circ^{q^2}},
\\\label{tipoF(4)-2}
& \xymatrix{ \circ^{-1} \ar@{-}[r]^{q} & \circ^{-1} \ar@{-}[r]^{q^{-2}} & \circ^{q^2} \ar@{-}[r]^{q^{-2}} & \circ^{q^2}},
\end{align}
\begin{align} \label{tipoF(4)-3}
& \xymatrix{ \circ^{-1} \ar@{-}[r]^{q^2} \ar@{-}[rd]_{q^{-1}} & \circ^{-1} \ar@{-}[r]^{q^{-2}} & \circ^{q^2} \\ & \circ^{q} \ar@{-}[u]_{q^{-1}} &},
\\ \label{tipoF(4)-4}
& \xymatrix{ \circ^{q^2} \ar@{-}[r]^{q^{-2}} & \circ^{-1} \ar@{-}[r]^{q^2} & \ar@{-}[ld]^{q^{-3}}\circ^{-1} \\ & \circ^{-1} \ar@{-}[u]^{q} &},
\\ \label{tipoF(4)-5}
& \xymatrix{ \circ^{q^2} \ar@{-}[r]^{q^{-2}} & \circ^{q} \ar@{-}[r]^{q^{-1}} & \circ^{-1} \ar@{-}[r]^{q^3} & \circ^{q^{-3}}},
\\ \label{tipoF(4)-6}
& \xymatrix{ \circ^{q^2} \ar@{-}[r]^{q^{-2}} & \circ^{q^2} \ar@{-}[r]^{q^{-2}} & \circ^{-1} \ar@{-}[r]^{q^3} & \circ^{q^{-3}}}.
\end{align}

\bigbreak
\noindent \textbf{Type $G(3)$}, $q \in \ku^\times$, $q^2, q^3 \neq 1$:
\begin{align}\label{tipoG(3)}
& \xymatrix{ \circ^{-1} \ar@{-}[r]^{q^{-1}} & \circ^{q} \ar@{-}[r]^{q^{-3}} & \circ^{q^3} },
\\ & \xymatrix{ \circ^{-1} \ar@{-}[r]^{q} & \circ^{-1} \ar@{-}[r]^{q^{-3}} & \circ^{q^3} },
\\ & \xymatrix{ & \circ^{-1} \ar@{-}[rd]^{q^3} & \\ \circ^{q} \ar@{-}[ru]^{q^{-1}} \ar@{-}[rr]^{q^{-2}} &  & \circ^{-1} },
\\ & \xymatrix{ \circ^{-q^{-1}} \ar@{-}[r]^{q^2} & \circ^{-1} \ar@{-}[r]^{q^{-3}} & \circ^{q^3} }.
\end{align}
\end{theorem}
\pf
When the braiding is of type $A_\theta$ or $B_\theta$ it follows by \cite[Propositions 3.9, 3.10]{Ang-st}. The proof for the other cases is completely analogous, so we just show in detail the case $C_\theta$. A first remark is that the submatrix $(q_{ij})_{1 \leq i,j \leq \theta-1}$ is of type $A_{\theta-1}$, so it is standard and has a generalized Dynkin diagram as \eqref{tipoA}, and the submatrix $\left( \begin{array}{cc} q_{\theta-1, \theta-1} & q_{\theta-1, \theta} \\ q_{\theta, \theta-1} & q_{\theta\theta} \end{array} \right)$ is of type $B_2$.

Also, if there exists $1 \leq i \leq \theta-1$ such that $p(\alpha_i)=-1$, then the reflection $s_i$ changes the set of roots, and by Remark \ref{obs:simetria trenza 2}, $q_{ii}=-1 \neq q_{i,i-1}q_{i-1,i}, q_{i,i+1}q_{i+1,i}$.

If $p(\alpha_i)=1$ for any $i$, then the root system is of finite type, and it follows that the braiding is of Cartan type by \cite[Prop. 3.8]{Ang-st}. If not, we can assume that $p(\alpha_{\theta-1})=-1$ up to applying a suitable sequence of reflections $s_i$. Applying $s_{\theta-1}$, $\alpha_\theta$ becomes odd for the new parity function, and moreover the reflection $s_{\theta}$ changes the root system, so in the original braiding $q_{\theta\theta}\neq -1$ by Remark \ref{obs:simetria trenza 2}.

We can make also $p(\alpha_{\theta-1})=1$ up to applying some reflections. In such case, $a_{\theta-1, \theta}=-2$, $a_{\theta-1, \theta-2}=-1$ and
applying the reflection $s_{\theta-1}$ the vertices $\theta-2$, $\theta$ are not connected (i.e. $\alpha_{\theta-2}+\alpha_\theta$ is not a root). In
consequence, if $(\widetilde q_{ij})$ denotes the braiding matrix after applying the reflection $s_{\theta-1}$, we have
\begin{align*}
1  & = \widetilde q_{\theta-2, \theta} \widetilde q_{\theta, \theta-2} = \chi (\alpha_{\theta-2}+\alpha_{\theta-1}, \alpha_{\theta}+ 2 \alpha_{\theta-1}) \chi (\alpha_{\theta}+ 2 \alpha_{\theta-1}, \alpha_{\theta-2}+\alpha_{\theta-1})
\\ & = q_{\theta-1, \theta-1}^4 q_{\theta-2, \theta-1}^2 q_{\theta-1, \theta-2}^2 q_{\theta, \theta-1}q_{\theta-1, \theta}.
\end{align*}
If we assume $q_{\theta-1, \theta-1}^2 \neq q_{\theta \theta}$, by the possible values of a matrix of type $B_2$ in such conditions we obtain $\widetilde q_{\theta-2, \theta} \widetilde q_{\theta, \theta-2} = -1$, a contradiction. Therefore $q_{\theta-1, \theta-1}^2 = q_{\theta \theta}$, and the braiding has a generalized Dynkin diagram as in \eqref{tipoC}. By direct computation we can see that a braiding as before is of type $C_\theta$.
\epf

\medskip

\begin{obs}
Note that these diagrams correspond with the following ones in Heckenberger's classification \cite{He}:
\begin{itemize}
  \item row 7 for type $G(3)$, and rows 9, 10 and 11 for type $D(2,1;\alpha)$ in Table 2,
  \item row 9 for type $F(4)$ in Table 3,
  \item rows 1 and 2 for type $A_\theta$, rows 3, 4, 5 and 6 for type $B_\theta$, rows 7, 8, 9 and 10 for type $C_\theta$, $D_\theta$ in Table 4.
\end{itemize}
\end{obs}

\subsection{Presentation by generators and relations}\label{subsec:gen-rels}

In this subsection we will present the Nichols algebras with super root systems of type $A$, $B$, $C$, $D$ by generators and relations using the results in \cite{Ang-rel}.

Recall that $[x,y]_c=xy-\chi(\alpha, \beta)yx$ if $x,y \in \toba(V)$ are homogeneous of degree $\alpha, \beta \in \zt$, respectively. In particular,
$(\ad_c x_i)(y):=[x_i, y]_c$.

We will define the hyperletters associated to the root vectors of type $A$, $B$, $C$, $D$, following \cite[Corollary 2.17]{Ang-rel}.

First of all, $x_{u_{ii}}=x_{\alpha_i}=x_i$, and recursively,
\begin{equation}\label{hyperpalabra u}
x_{\ub_{ij}}:= [x_i, x_{\ub_{i+1,j}}]_c, \quad i<j.
\end{equation}
Also, $x_{\vb_{i,\theta}}= [x_{\ub_{i,\theta}}, x_{\theta}]_c$, and recursively,
\begin{equation}\label{hyperpalabra v}
x_{\vb_{ij}}:= [x_{\vb_{i,j+1}}, x_j]_c, \quad i<j.
\end{equation}
For type $C$, $x_{\wb_{i,\theta-1}}= [x_{\ub_{i,\theta}}, x_{\theta-1}]_c$, and then,
\begin{align}\label{hyperpalabra w}
x_{\wb_{ij}}:= [x_{\wb_{i,j+1}}, x_j]_c, \quad i<j,
\\ x_{\widetilde \wb_i}:= [x_{\ub_{i,\theta-1}}, x_{\ub_{i,\theta}}]_c\label{hyperpalabra w2}.
\end{align}
For type $D$, we have $x_{\widetilde \ub_{\theta-2}}=[x_{\theta-2}, x_\theta]$, and recursively $x_{\widetilde \ub_i}=[x_i, x_{\widetilde \ub_{i+1}}]$. Also, $x_{\zb_{i,\theta-2}}= [x_{\ub_{i,\theta}}, x_{\theta-2}]_c$, and
\begin{align}\label{hyperpalabra z}
x_{\zb_{ij}}:= [x_{\zb_{i,j+1}}, x_j]_c, \quad i<j,
\\ x_{\widetilde \zb_i}:= [x_{\ub_{i,\theta-1}}, x_{\widetilde \ub_{i}}]_c\label{hyperpalabra z2}.
\end{align}
In this case, note that $x_{\ub_{\theta-2, \theta}}= [[x_{\theta-2}, x_\theta]_c, x_{\theta-1}]_c$.
\medskip

For any $\alpha \in \de^V_+$ we write $q_{\alpha}= \chi(\alpha, \alpha)$, and $N_\alpha= \ord (q_\alpha)$.

\begin{theorem}\label{theorem:presentation} Let $(V,c)$ be a braided vector space of diagonal type, with braiding matrix $(q_{ij})_{1 \leq i,j \leq \theta}$. Assume that the root system of $\toba (V)$ is of super type, with connected components of type $A$, $B$, $C$, $D$.

\noindent The Nichols algebra $\toba(V)$ is presented by generators $x_i$, $1 \leq i \leq \theta$, and relations:
\begin{align}
& x_{\alpha}^{N_\alpha}=0, \quad \alpha \in \de^V_+ \label{powerrootvector}, \\
& (\ad_c x_i)^{1-a_{ij}} x_j=0, \quad q_{ii}^{1-a_{ij}} \neq 1; \label{quantum Serre}
\end{align}
if $q_{kk}=-1$, $a_{kj}=a_{kl}=-1$ and $q_{kj}q_{jk}q_{kl}q_{lk}=q_{jl}q_{lj}=1$, then
\begin{equation}\label{relation A}
[\ad_c x_j \ad_c x_k (x_l), x_k]_c=0;
\end{equation}
if  there exist components of type $B_N$ such that $q_{NN}^2+q_{NN}+1=q_{N-1,N-1}+1=0$, then
\begin{align}\label{relation B1}
& [x_{\ub_{N-1,N}}, \  x_{\vb_{N-1,N}}]_c=0;
\\ \label{relation B2}
& [x_{\vb_{N-2,N}}, \  x_{\ub_{N-1,N}}]_cc=0;
\end{align}
if there exist components of type $C_N$, then
\begin{align}
& [x_{\widetilde \wb_{N-2}}, x_{N-1}]_c=0, \qquad q_{N-2, N-2}=q_{N-1, N-1}=-1, \label{relation C1}
\\ & [x_{\wb_{N-3,N-2}}, x_{N-1}]_c=0, \qquad q_{N-2, N-2} \neq q_{N-1, N-1}=-1, \label{relation C2}
\\ & [x_{\wb_{N-2, N-1}}, x_{N-1}]_c=0, \qquad q_{N-1, N-1}\in \G_3; \label{relation C3}
\end{align}
if there exist components of type $D_N$ and $q_{N-1,N-1}=q_{NN}=-1$, then
\begin{equation}\label{relation D}
[[x_{N-2}, x_{N-1}]_c, x_N]_c+  q_{N-2,N-1}q_{N-1,N-2} q_{N-1,N} \  x_{\ub_{N-2, N}}=0.
\end{equation}
\end{theorem}
\pf It is enough to consider the connected case. In such case, \cite[Theorem 3.9]{Ang-rel} gives us a family of relations which generates the ideal of relations. This set contains redundant relations, so we work as in \cite[Sections 5B, 5C]{Ang-st} to prove that the relations in large rank are generated by the relations in small rank analyzing each case.
\epf

\begin{obs}\label{obs:redundant power root vectors}
Some of the relations \eqref{powerrootvector} are redundant. Depending on the type we can restrict the set of such relations to the following sets:
\begin{align}
x_\alpha^{N_\alpha} &=0, \qquad \alpha\in\Delta_+^V, \ \chi(\alpha,\alpha)\ne -1,
\\ x_i^2 &=0, \qquad i\in I,\ q_{ii}\ne -1,
\end{align}
for diagrams \eqref{tipoA}, \eqref{tipoB2}, \eqref{tipoB3}, \eqref{tipoC}, \eqref{tipoD1}, \eqref{tipoD2};
\begin{align}
x_2^3 = x_\alpha^{N_\alpha} =0, \qquad \alpha\in\{\alpha_1, \alpha_1+2\alpha_2\},
\end{align}
for diagram \eqref{tipoB1}.
\end{obs}

\begin{obs}\label{obs:relaciones Yamane}
Notice the similarity of this presentation to the one for quantized enveloping superalgebras in \cite{Y}. Anyway, when $q$ has small order, we need some extra relations. We obtain also a presentation for non-symmetric matrices $(q_{ij})$.
\end{obs}

\begin{obs}\label{obs:presentacion D(2,1;alpha)}
The presentation by generators and relations of Nichols algebras of type $D(2,1;\alpha)$ is given in \cite[Propositions 4.1, 4.2]{Ang-rel}.
\end{obs}

\begin{obs}\label{obs:gen-deg-one-new} As said in Remark \ref{obs:gen-deg-one},
the Conjecture on generation in degree one was verified in many cases.
The pattern of all the proofs is the same \cite{AS-ann-ec-sci}: to analyze the braided subspaces of the tensor algebra spanned by one of the relations and the
elements of the basis intervening on it; these braided subspaces turns out to have an infinite-dimensional Nichols algebra by \cite{He}.
Here again, it can be shown that \fd{} pointed Hopf algebras with infinitesimal braiding of super type are generated in degree one. The proof
is analogous to that of \cite[Th. 2.7]{Ang-GI}. It follows that a pointed Hopf algebra whose infinitesimal braiding has connected components of rank greater or equal than 8 is generated in degree one.
\end{obs}

\end{document}